\documentclass[12pt]{article}

\usepackage{amstext    }
\usepackage{amsthm    }
\usepackage{a4}
\usepackage[mathscr]{eucal}
\usepackage{mathrsfs}

\usepackage{amsmath}
\usepackage{amssymb}
\usepackage{amscd}
\newtheorem{theorem}{Theorem}[section]
\newtheorem{definition}[theorem]{Definition}
\newtheorem{proposition}[theorem]{Proposition}
\newtheorem{corollary}[theorem]{Corollary}
\newtheorem{lemma}[theorem]{Lemma}

\newtheorem{example}[theorem]{Example}

\newtheorem{problem}[theorem]{Problem}

\newcommand{\cali}[1]{\mathscr{#1}}

\newcommand{\const}{\mathop{\mathrm{const}}}

\newcommand{\dist}{\mathop{\mathrm{dist}}\nolimits}

\newcommand{\ddc}{dd^c}

\newcommand{\dbar}{\overline\partial}

\newcommand{\capacity}{\mathop{\mathrm{cap}}\nolimits}

\newcommand{\BTK}{{\rm BTK}}

\newcommand{\Cc}{\cali{C}}

\newcommand{\Uc}{\cali{U}}
\newcommand{\Vc}{\cali{V}}

\newcommand{\Tc}{\cali{T}}

\newcommand{\B}{\mathbb{B}}
\newcommand{\C}{\mathbb{C}}
\newcommand{\N}{\mathbb{N}}

\newcommand{\R}{\mathbb{R}}
\renewcommand\P{\mathbb{P}}

\title{Characterization of Monge-Amp\`ere measures with H\"older continuous potentials}

\author{Tien-Cuong Dinh and Vi{\^e}t-Anh Nguy{\^e}n}

\begin{document}

\maketitle

\begin{abstract}
We show that the complex Monge-Amp\`ere equation on a compact K\"ahler manifold $(X,\omega)$ of dimension $n$ admits a H\"older continuous $\omega$-psh solution if and only if its right-hand side is a positive measure with H\"older continuous super-potential. This property is true in particular when the measure has locally H\"older continuous potentials or when it belongs to the Sobolev space $W^{2n/p-2+\epsilon,p}(X)$ or to the Besov space $B^{\epsilon-2}_{\infty,\infty}(X)$ for some $\epsilon>0$ and $p>1$. 
\end{abstract}

\noindent
{\bf Classification AMS 2010}: 32U, 32W20, 53C55.


\section{Introduction} \label{introduction}

 Let $(X,\omega)$ be a compact K\"ahler manifold of dimension $n$.  Recall that   a function
$u:\ X\to \R\cup\{-\infty\}$ is  said  to be {\it  $\omega$-psh} if it is locally the difference of a psh function and a potential of $\omega$, i.e. a smooth function $v$ such that $\ddc v=\omega$. We consider the complex Monge-Amp\`ere equation 
$$(\ddc\varphi+\omega)^n=\mu,$$
where $\mu$ is a positive measure on $X$ and $\varphi$ a bounded $\omega$-psh function, see Bedford-Taylor \cite{BedfordTaylor}, Demailly \cite{Demailly} and Forn\ae ss-Sibony \cite{FornaessSibony} for the intersection of currents and for basic properties of psh functions. For cohomology reason, the above relation implies that the mass of $\mu$ is equal to the mass of the measure $\omega^n$, i.e.
$$\|\mu\|=\int_X\omega^n.$$
In what follows, we always assume this condition.

When $\mu$ is a smooth volume form, the famous theorem of Calabi-Yau says that the Monge-Amp\`ere equation admits a smooth solution $\varphi$ which is unique up to an additive constant \cite{Calabi,Yau}. In this paper, we consider the case with H\"older continuous solutions. 
Without reviewing the long history of the complex Monge-Amp\`ere equation,
let us mention few steps in the recent development.

In \cite{Kolodziej1,Kolodziej2}, Kolodziej has constructed a continuous solution under some hypothesis on the measure $\mu$, in particular for $\mu$ of class $L^p$, $p>1$. 
Then, he proved in \cite{Kolodziej3} that the solution is H\"older continuous when $\mu$ is of class $L^p$, $p>1$, see also \cite{EGZ,GuedjKolodziejZeriahi}. 
Some important steps in his approach were improved by Dinew and Zhang \cite{DinewZhang}.
Very recently based on Demailly's regularization method \cite{Demailly2,Demailly3} and the above Dinew-Zhang's results, Demailly, Dinew, Guedj, Hiep, Kolodziej, Zeriahi \cite{DemaillyDinew} obtained an explicit H\"older exponent of the solution, see also \cite{Dinew}. A necessary condition on $\mu$ to have a H\"older continuous solution $u$ was obtained by Dinh-Nguyen-Sibony in \cite{DinhNguyenSibony}.
We refer to the works of the above mentioned authors and Eyssidieux, Pali, Plis, Song, Tian \cite{DemaillyDinew, DemaillyPali, Dinew, DinewZhang, DinhNguyenSibony,EGZ,Hiep,Kolodziej1,Kolodziej2,Kolodziej3,Plis, SongTian,TianZhang} for results in this direction, for related topics and a more complete list of references.

Here is our main theorem which implies several known results. The proof uses the above results and some ideas from the works by Sibony and the authors  \cite{DinhSibony1, DinhSibony2, DinhSibony3,DinhSibony4}.

\begin{theorem} \label{th_main}
Let $(X,\omega)$ be a compact K\"ahler manifold of dimension $n$. Let $\mu$ be a positive measure on $X$ of mass $\int_X\omega^n$. Then the Monge-Amp\`ere equation
$$(\ddc\varphi+\omega)^n=\mu$$
admits a H\"older continuous $\omega$-psh solution $\varphi$ if and only if $\mu$ admits a H\"older continuous super-potential. 
\end{theorem}

The H\"older exponent of the solution depends on $n$ and on the H\"older exponent of the super-potential $\Uc$ of $\mu$. It will be specified in the proof of Theorem \ref{th_main} in Section \ref{section_sp}. 

We will recall the notion of super-potentials for measures in Section \ref{section_sp}. Super-potentials for positive closed currents were introduced by Sibony and the first author in \cite{DinhSibony2,DinhSibony3, DinhSibony4}. They play an important role in complex dynamics. In the references above, situations where the H\"older continuity of super-potential can be verified, are described. The reader will also find in Section \ref{section_sp} some methods to obtain this property. In particular, we will see that the class of measures with
H\"older continuous super-potential contains the measures given by $L^p$ forms, $p>1$, considered in \cite{DemaillyDinew, DinewZhang,Kolodziej3} and the measures satisfying some Hausdorff type regularity investigated  in \cite{Hiep}. 

The following corollary is a consequence of Theorem \ref{th_main} and Proposition  \ref{prop_example} below.

\begin{corollary} \label{main_cor}
Under the hypotheses of Theorem \ref{th_main}, if we can write locally $\mu=\ddc U+\partial V+\dbar W$ with H\"older continuous forms $U,V,W$ of bidegrees $(n-1,n-1)$, $(n-1,n)$ and $(n,n-1)$ respectively, then the considered Monge-Amp\`ere equation admits a H\"older continuous $\omega$-psh solution.  Moreover, the hypothesis on $\mu$ is satisfied when this measure belongs to the Sobolev space $W^{2n/p-2+\epsilon,p}(X)$ or to the Besov space $B^{\epsilon-2}_{\infty,\infty}(X)$  for some $\epsilon>0$ and $p>1$. 
\end{corollary}

In the case with parameters $(X_t,\omega_t,\mu_t)$, where the compact K\"ahler manifolds $(X_t,\omega_t)$ have uniformly bounded geometry and the super-potentials of the measures $\mu_t$ are uniformly H\"older continuous, we obtain solutions $\varphi_t$ which are uniformly H\"older continuous. We can also extend our results to the case of a big and nef class with a solution locally H\"older continuous on the ample locus. We refer to \cite{DemaillyDinew,DinewZhang, EGZ} for the techniques needed to handle these situations.

The following problem suggested by the work of Sibony and the authors \cite{DinhNguyenSibony} is still open. We refer to \cite{DemaillyDinew,DinhSibony1} for some particular cases where the answer is positive.

\begin{problem}
Let $\mu$ be a probability measure on $X$. Assume that $\mu$ is moderate. Does the Monge-Amp\`ere equation 
$$(\ddc\varphi+\omega)^n=\mu$$
admit a H\"older continuous $\omega$-psh solution $\varphi$ ?  
\end{problem}

The notion of moderate measures will be recalled in Section \ref{section_moderate}. We think that the answer is not always positive and the problem requires probably a better understanding of the notion of capacity $\Tc(\cdot)$, see \cite{DinhSibony5} and Section   \ref{section_moderate}  for the definition.

\begin{problem} 
Characterize the positive measures $\mu$ on $X$ such that the associated complex
Monge-Amp\`ere equation admits a continuous (resp. bounded) solution.
\end{problem}

Note that it is not difficult to show using \cite[Ch.III (3.6) and (3.11)]{Demailly}
that when $\varphi$ is a continuous (resp. bounded) $\omega$-psh function the measure $(\ddc \varphi+\omega)^n$ has a continuous (resp. bounded) super-potential. However, in dimension $n\geq 2$, the last property does not characterize the Monge-Amp\`ere measures with continuous (resp. bounded) potential. We can consider, for example in dimension 2, a measure with a single singularity likes $(\ddc \log(-\log\|z\|))^2$. It has a continuous super-potential.

\bigskip
\noindent
{\bf Acknowledgement.} The authors thank Nicolas Lerner and Nessim Sibony for their help during the preparation of this work.

\section{Moderate measures and capacities} \label{section_moderate}

In what follows,  $(X,\omega)$ always denotes  a  compact K\"{a}hler manifold of dimension $n$.
Recall from  \cite{DinhNguyenSibony,DinhSibony1} that a positive measure $\mu$  on $X$  is {\it moderate} if there are constants $c>0$ and $\alpha>0$ such that 
$$\int e^{-\alpha u}d\mu\leq c \quad \mbox{for every } u \ \omega \mbox{-psh such that } \int_X u\omega^n=0.$$
Note that the functions $u$ satisfying the last condition describe a compact set of $\omega$-psh functions. The condition can be replaced by other ones, e.g. $\max_X u=0$.

\begin{lemma} \label{lemma_sublevel}
Let $\mu$ be a moderate positive measure on $X.$ Then there are constants $c>0$ and $\alpha>0$ such that for any $M\geq 0$ 
$$\mu\{u<-M\} \leq ce^{-\alpha M} \quad \mbox{for every } \omega \mbox{-psh function } u \mbox{ such that } \int_X u\omega^n=0.$$
Moreover, if $p\geq 1$ is a real number, then there is a constant $c_p>0$  such that 
$$\|u\|_{L^p(\mu)}\leq c_p \quad \mbox{and} \quad \int_{\{u<-M\}} |u|^pd\mu \leq c_pe^{-\alpha M} \quad \mbox{for all } M \mbox{ and } u \mbox{ as above}.$$
\end{lemma}
\proof
The first assertion follows from the definition of moderate measure. Observe that $u$ is bounded above by a constant independent of $u$. Therefore, the second assertion in the lemma follows from the first one and  the identities
$$    \int_{|u|>M}  (|u|^p-M^p)d\mu=\int_{M^p}^\infty \mu\{ |u|^p>r\}dr=p\int_M^\infty \mu\{ |u|>t\}t^{p-1}dt$$
for $M\geq 0$.  We change $\alpha$ if necessary.
\endproof

\begin{lemma} \label{lemma_moderate_Lp}
Let $\mu$ be a moderate positive measure on $X$. Then for any real number $p\geq 1$, there is a constant $c>0$ such that 
$$\|u-u'\|_{L^p(\mu)}\leq c\max(1,-\log\|u-u'\|_{L^1(\mu)})^{(p-1)/p} \|u-u'\|_{L^1(\mu)}^{1/p}$$
for all $\omega$-psh functions $u$ and $u'$ such that $\int_Xu\omega^n=\int_Xu'\omega^n=0$. 
\end{lemma}
\proof
In what follows,  $\lesssim$ and $\gtrsim$ denote inequalities up to a multiplicative constant.
Since  $\|u\|_{L^p(\mu)}$ and $\|u'\|_{L^p(\mu)}$ are bounded, we only have to consider the case where $\|u-u'\|_{L^1(\mu)}$ is small. Define for a constant $M$ large enough $u_M:=\max(u,-M)$ and $u_M':=\max(u',-M)$. Using Lemma \ref{lemma_sublevel} and H\"older's inequality, we have 
\begin{eqnarray*}
\|u-u'\|_{L^p(\mu)} & \lesssim & \|u_M-u'_M\|_{L^p(\mu)} + e^{-\alpha M}\lesssim M^{(p-1)/p}\|u_M-u_M'\|^{1/p}_{L^1(\mu)}+e^{-\alpha M} \\
& \lesssim & M^{(p-1)/p}\|u-u'\|^{1/p}_{L^1(\mu)}+M^{(p-1)/p}e^{-\alpha M/p}+e^{-\alpha M}.
\end{eqnarray*}
It is enough to choose $M$ equal to a large constant times $-\log\|u-u'\|_{L^1(\mu)}$. 
\endproof

Recall the following notion of capacity which was introduced by Kolodziej \cite{Kolodziej2} and is related to the well-known Bedford-Taylor capacity \cite{BedfordTaylor}. For any Borel subset $A$ of $X$, define 
$$\capacity_\BTK(A):=\sup \Big\{\int_A(\ddc u+\omega)^n, \quad u \ \omega \mbox{-psh such that } 0\leq u\leq 1\Big\}.$$
The following definition is inspired by the works of Kolodziej \cite{Kolodziej1, Kolodziej2, Kolodziej3}.

\begin{definition} \rm \label{def_K_moderate}
A positive measure $\mu$ on $X$ is said to be {\it K-moderate} if there are constants $c>0$ and $\alpha>0$ such that for every Borel subset $A$ of $X$ we have
$$\mu(A)\leq c \exp(-\capacity_\BTK(A)^{-\alpha}).$$
\end{definition}

We will also need the following notion of capacity introduced by Sibony and the first author \cite{DinhSibony5} which is related to the capacities of Alexander  \cite{Alexander} and of Siciak  \cite{Siciak}, see also Harvey-Lawson \cite{HarveyLawson}.  For any Borel subset $A$ of $X$, define 
$$\Tc(A):=\inf \Big\{\exp(\sup_A u),\quad u \ \omega \mbox{-psh and } \max_Xu=0\Big\}.$$ 
 Recall   from Guedj-Zeriahi \cite[Prop. 7.1]{GuedjZeriahi} 
the  following  relation  between  $\capacity_\BTK(A)$ and  $\Tc(A)$  for every compact set $A\subset X:$  
\begin{equation}\label{eq_capacities}
c_1\exp(-\lambda_1 \capacity_\BTK(A)^{-1})\leq \Tc(A)\leq c_2\exp(-\lambda_2\capacity_\BTK(A)^{-1/n}),
\end{equation}
where $c_i>0$ and $\lambda_i>0$ are constants independent of $A.$ 

\begin{proposition} \label{prop_moderate_K}
Let $\mu$ be a probability measure on $X$. Then $\mu$ is K-moderate if and only if it is weakly moderate, i.e., there are constants  $\lambda>0$ and $\alpha>0$ such that
$$\int \exp(\lambda |u|^\alpha) d\mu \leq c \quad \mbox{for } u \ \omega\mbox{-psh with } \max_X  u=0.$$
In particular, if $\mu$ is moderate, then it is K-moderate.
\end{proposition}
\proof
Assume that $\mu$ is weakly moderate as above. We  will show that it is K-moderate. 
It is enough to obtain the estimate in Definition \ref{def_K_moderate} for every compact set $A$  such that  $\capacity_\BTK(A)$ is  small enough. 
We can also assume that $\capacity_\BTK(A)$ is strictly positive since otherwise $A$ is pluripolar and the fact that $\mu$ is weakly moderate implies that $\mu(A)=0$. So we also have $\Tc(A)>0$. 

Let $u$ be an $\omega$-psh function with $\max_X u=0$ such that 
$$u\leq \log\Tc(A)+1 \mbox{ on } A.$$
Using that $\Tc(A)$ is small, we obtain that
$$\mu(A)\leq \mu\big\{u\leq \log\Tc(A)+1\big\} \leq \mu\big\{\lambda|u|^\alpha\geq (-\log \Tc(A))^{\alpha/2}\big\}.$$
Since $\mu$ is weakly moderate, if $c$ is the constant as in the lemma, the last quantity is bounded above by 
$$c\exp\big(-(-\log\Tc(A))^{\alpha/2}\big), $$
which is, by the   second inequality in (\ref{eq_capacities}),  dominated  by
 $$ c'\exp\big(-\capacity_\BTK(A)^{-\alpha'}\big)$$
for some constants $c'>0$ and $\alpha'>0$. So $\mu$ is K-moderate.

Assume now that $\mu$ is K-moderate as in Definition \ref{def_K_moderate}. We  will show that it is weakly moderate. Let $u$ be an $\omega$-psh function such that $\max_X u=0$. It is enough to show for any $M\geq 0$ that 
$$\mu\{u< -M\} \leq c'' \exp(-M^{\alpha''})$$
for some constants $c''>0$ and $\alpha''>0$ independent of $M$ and of $u$. Let $A$ be an arbitrary compact subset of the open set $\{u<-M\}$. We want to bound $\mu(A)$. We have $\Tc(A)\leq e^{-M}$ by definition of $\Tc$. 
We can assume that $\capacity_\BTK(A)$ and $\Tc(A)$ are small enough.
Since $\mu$ is K-moderate, we obtain for some constant  $\alpha''>0$
$$\log\mu(A)\leq \const -\capacity_\BTK(A)^{-\alpha}\leq -\big(-\log\Tc(A)\big)^{\alpha''} \leq -M^{\alpha''},$$
where the  second  inequality  follows  from the first   estimate  in  (\ref{eq_capacities}).
Hence, $\mu$ is weakly moderate. The lemma follows.
\endproof

\begin{example} \rm
Let $\mu$ be a probability measure on $\P^1$ smooth except at 0 and such that 
$$\mu=|z|^{-2}\exp\big(-(-\log|z|)^{1/2}\big)(idz\wedge d\overline z) \quad \mbox{near }0.$$ 
If $u$ is equal to $\log|z|$ near 0, we see that $\exp(-\lambda u)$ is not $\mu$-integrable for any $\lambda>0$. So $\mu$ is not moderate. One can check that $\mu$ is, however, weakly moderate and hence it is K-moderate. 
\end{example}

\section{Super-potentials of positive measures} \label{section_sp}

The notion of super-potential was introduced by Sibony and the first author. It extends the notion of potential of positive closed $(1,1)$-currents to positive closed $(p,p)$-currents and allows to solve some problems in complex dynamics. We recall it in the case of measures, i.e. for $p=k$. 

Let $\Cc$ denote the set of positive closed $(1,1)$-currents in the cohomology class $\{\omega\}$. This is a convex compact set. We can consider for each real number $\alpha>0$ the following distance on $\Cc$:
$$\dist_\alpha(T,T'):=\sup_{\|\Phi\|_{\Cc^\alpha}\leq 1} |\langle T-T',\Phi\rangle|$$
where $\Phi$ is a test smooth $(n-1,n-1)$-form. Observe that the family $\dist_\alpha$ is decreasing in $\alpha$.
The following proposition was obtained in \cite{DinhSibony2} as a consequence of the interpolation theory for Banach spaces. 

\begin{proposition} \label{prop_norm_alpha}
Let $\alpha$ and $\beta$ be real numbers such that $\beta\geq\alpha>0$. Then there is a constant $c=c(\alpha,\beta)>0$ such that 
$$\dist_{\beta}\leq \dist_\alpha \leq c(\dist_\beta)^{\alpha/\beta}.$$
\end{proposition}

For each current $T$ in $\Cc$ there is a unique $\omega$-psh function $u$ such that
$$T=\ddc u +\omega \quad \mbox{and} \quad \int_X u \omega^n=0.$$
We call $u$ the {\it $\omega$-potential} of $T$.
For each real number $p>1$, define the following distance on $\Cc$
$$\dist_{L^p}(T,T'):=\|u-u'\|_{L^p},$$
where $u$ and $u'$ are $\omega$-potentials of $T$ and $T'$ respectively and the $L^p$ norm is with respect to the measure $\omega^n$. Since we assume that $\omega^n$ is a probability measure, the 
family $\dist_{L^p}$ is increasing in $p$. Note that $\Cc$ has finite diameter with respect to all the above distances.

\begin{proposition} \label{prop_norm_lp}
Let $p\geq 1$ be a real number. Then there are constants $c>0$, $c'>0$, and $c''>0$ depending on $p$ such that 
$$c\dist_2 \leq \dist_{L^1} \leq c' \dist_1 \quad \mbox{and} \quad \dist_{L^p}\leq c''\max(1,-\log\dist_{L^1})^{{p-1\over p}}(\dist_{L^1})^{1\over p}.$$ 
\end{proposition}
\proof
Given $T,T'\in\Cc,$   let $u$ and $u'$ be the $\omega$-potentials of $T$ and $T'$ respectively.
Let $\Phi$ be a smooth $(n-1,n-1)$-form such that $\|\Phi\|_{\Cc^2}\leq 1$. We have
$$|\langle T-T',\Phi\rangle|=|\langle u-u',\ddc\Phi\rangle|\lesssim \dist_{L^1}(T,T').$$
So the first inequality in the proposition is clear.

Consider the second inequality in the proposition. Let $\pi_1$ and $\pi_2$ denote the canonical projections from $X\times X$ onto its factors. Let $\Delta$ be the diagonal of $X\times X$. 
In the proof of \cite[Prop. 2.1]{DinhSibony6}, an explicit kernel $K(x,y)$ on $X\times X$ was obtained, see also Bost-Gillet-Soul\'e \cite{BostGilletSoule}. It is an $(n-1,n-1)$-form smooth outside $\Delta$ such that 
$$\|K(x,y)\|\lesssim -\log\dist(x,y)\dist(x,y)^{2-2n}$$ 
and
$$\|\nabla K(x,y)\|  \lesssim  \dist(x,y)^{1-2n}$$
when $(x,y)$ tends to $\Delta$. 
Here, $\|\nabla K(x,y)\|$ denotes the sum of the norms of the gradients of the coefficients of $K(x,y)$ for a fixed atlas on $X\times X$. 

This kernel gives a solution $v$ to the equation $\ddc v=T-T'$ in the sense of currents with
$$v(x):=\int_{y\in X} K(x,y)\wedge (T(y)-T'(y))$$ 
or more formally
$$v:=(\pi_1)_*(K\wedge \pi_2^*(T-T')).$$
Indeed, if $[\Delta]$ denotes the current of integration on the diagonal $\Delta$, then $[\Delta]-\ddc K$ is a smooth representation of the cohomology class of $\Delta$ in the K\"unneth decomposition of $H^{n,n}(X\times X,\C)$.

We  show  that   $\|v\|_{L^1}\lesssim \dist_1(T,T')$. Consider a test 
smooth $(n,n)$-form $\Phi$ such that $\|\Phi\|_\infty\leq 1$. We have
 $$\langle v,\Phi\rangle =\int_{X\times X} K(x,y)\wedge\Phi(x)\wedge (T(y)-T'(y))= \big\langle T-T',(\pi_2)_*(K\wedge \pi_1^*(\Phi))\big\rangle.$$
The estimates on $K$ imply that  $(\pi_2)_*(K\wedge \pi_1^*(\Phi))$ is a form with bounded $\Cc^1$ norm. 
We deduce that $|\langle v,\Phi\rangle|\lesssim \dist_1(T,T')$. Since this property holds for every $\Phi$, we obtain that  $\|v\|_{L^1}\lesssim \dist_1(T,T')$.

Define $m:=\int_X v \omega^n$ and $\widetilde v:=v-m$. It follows  from the above  discussion  that  $|m|\lesssim \dist_1(T,T')$ and $\ddc \widetilde v=T-T'=\ddc (u-u')$. Since the solution of the equation $\ddc \widetilde v=T-T'$ with $\int_X\widetilde v\omega^n=0$ is unique, we obtain that 
$\widetilde v=u-u'$.  This, combined  with  the  estimates $\|v\|_{L^1}\lesssim\dist_1(T,T')$ and  $|m|\lesssim \dist_1(T,T')$,
implies that $\|u-u'\|_{L^1}\lesssim \dist_1(T,T')$. The second inequality in the proposition follows.

Recall that the measure $\omega^n$ is moderate, see \cite{Skoda,Zeriahi}. So Lemma \ref{lemma_moderate_Lp} applied to $\omega^n$ implies the last inequality in the proposition. 
\endproof

Propositions \ref{prop_norm_alpha} and \ref{prop_norm_lp} show that the distances considered above define the same topology on $\Cc$. It is not difficult to see that this topology is induced by the weak topology on currents. The above propositions also imply that
the following notion does not depend on the choice of the distance on $\Cc$ and therefore gives us a large flexibility to prove this H\"older property.

\begin{definition}\rm
A function $\Uc:\Cc\to \R$ is said to be {\it H\"older continuous} if it is H\"older continuous for one of the above distances on $\Cc$.
\end{definition}

\begin{definition}\rm
Let $\mu$ be a positive measure on $X$. {\it The super-potential} of $\mu$ is the function $\Uc:\Cc\to\R\cup\{-\infty\}$ defined by
$$\Uc(T):=\int u d\mu$$
where $u$ is the $\omega$-potential of $T$. 
\end{definition}

The infinite dimensional compact space $\Cc$ admits some ``complex structure" and
the super-potential $\Uc$ satisfies similar properties as the ones of quasi-psh functions in the finite dimensional case. However, we do not need these properties here. It is important to notice that the H\"older continuity, the continuity and the boundedness of the super-potential $\Uc$ do not depend on the choice of the K\"ahler form $\omega$ on $X$.

In what follows  we study  various  sufficient conditions  for   a  positive measure $\mu$ to possess  a H\"older continuous super-potential. We will need the following elementary lemma.

\begin{lemma} \label{lemma_elem}
Let $\mu$ be a positive measure on $X$. Then its super-potential $\Uc$ is H\"older continuous with H\"older exponent $0<\beta\leq 1$ with respect to the distance $\dist_{L^1}$ on $\Cc$ if and only if there is a constant $c>0$ such that 
$$\|u-u'\|_{L^1(\mu)}\leq c\max\big(\|u-u'\|_{L^1},\|u-u'\|_{L^1}^\beta\big)$$
for all $\omega$-psh functions $u$ and $u'$.
\end{lemma}
\proof
We first prove the necessary condition. Assume that $\Uc$ is H\"older continuous with exponent $\beta$ as above. 
Observe that it is enough to show that
$$\Big|\int (u-u')d\mu\Big|\lesssim  \max\big(\|u-u'\|_{L^1},\|u-u'\|_{L^1}^\beta\big).$$
Indeed, this inequality applied to $u,\max(u,u')$ and then to $u',\max(u,u')$ gives the result.

Consider first the case where $\int_Xu\omega^n=\int_Xu'\omega^n=0$. Define $T:=\ddc u+\omega$ and $T':=\ddc u'+\omega$. Then, by hypothesis, we have
$$\Big|\int (u-u')d\mu\Big|=|\Uc(T)-\Uc(T')|\lesssim \dist_{L^1}(T,T')^\beta=\|u-u'\|_{L^1}^\beta.$$
In the general case, define $m:=\int_Xu\omega^n$ and $ m':=\int_X u'\omega^n$. We can apply the first case to $v:=u-m$ and $v':=u'-m'$. In order to obtain the inequality in the lemma, it is enough to use the triangle inequality and to observe that $|m-m'|\lesssim \|u-u'\|_{L^1}$. 

For the sufficient part, assume the inequality in the lemma. Consider two currents $T$ and $T'$ in $\Cc$. Denote by $u$ and $u'$ their $\omega$-potentials which belong to a fixed compact family of $\omega$-psh functions. Therefore, $\|u\|_{L^1}$, $\|u'\|_{L^1}$ and $\|u-u'\|_{L^1}$ are bounded. The inequality in the lemma implies that 
$$|\Uc(T)-\Uc(T')|=\Big|\int (u-u')d\mu\Big| \lesssim \|u-u'\|_{L^1}^\beta.$$
The lemma follows.
\endproof

\begin{proposition} \label{prop_necessary}
Let $\varphi$ be a H\"older continuous $\omega$-psh function on $X$. Then the measure $\mu:=(\ddc\varphi+\omega)^n$ has a H\"older continuous super-potential.
\end{proposition}
\proof
Define $\mu_k:=(\ddc\varphi+\omega)^k\wedge \omega^{n-k}$. We prove by induction on $k$ that $\mu_k$ has a H\"older continuous potential.  Assume this is true for $k-1$. 
We will use the criterium given in Lemma \ref{lemma_elem}.
We can assume $u\geq u'$ since we can always reduce the problem to
the case with $u,\max(u,u')$ and the case with $u',\max(u,u')$. Subtracting from $u$ and $u'$ a constant allows to assume that $\int_X u\omega^n=0$. 

Define also $m':=\int_X u'\omega^n$ and $\widehat u':=u'-m'$. So $u$ and $\widehat u'$ belong to a fixed compact family of $\omega$-psh functions. We deduce that $\|u\|_{L^1}$ and $\|\widehat u'\|_{L^1}$ are bounded. 
By  \cite[Ch.III (3.11)]{Demailly}, the integrals $\langle\mu,u\rangle$ and  $\langle\mu,\widehat u'\rangle$ are also bounded. Therefore, by Lemma \ref{lemma_elem}, we only have to consider the case where $|m'|$ is bounded by a fixed constant large enough and to prove the inequality
 $$\int(u-u')d\mu_k\lesssim \|u-u'\|_{L^1}^{\beta_k}$$
for some constant $\beta_k>0$. Indeed, $\|u-u'\|_{L^1}$ is bounded by a fixed constant.

Since $\varphi$ is H\"older continuous, using a standard convolution and a partition of unity, we can write  $\varphi=\varphi_\epsilon+(\varphi-\varphi_\epsilon)$ with $\|\varphi_\epsilon\|_{\Cc^2}\lesssim \epsilon^{-2}$ and $|\varphi-\varphi_\epsilon|\lesssim \epsilon^\alpha$ for some $\alpha>0$. We have for $T:=\ddc u+\omega$ and $T':=\ddc u'+\omega$
\begin{eqnarray*}
|\langle \mu_k,u-u'\rangle|
 & \leq &  |\langle \mu_{k-1},u-u'\rangle|
 + |\langle \ddc\varphi\wedge (\ddc\varphi+\omega)^{k-1}\wedge \omega^{n-k},u-u'\rangle|\\
& \leq &   |\langle \mu_{k-1},u-u'\rangle|
+  |\langle \ddc\varphi_\epsilon\wedge (\ddc\varphi+\omega)^{k-1}\wedge \omega^{n-k},u-u'\rangle|\\
& & +  |\langle (\varphi-\varphi_\epsilon)\wedge (\ddc\varphi+\omega)^{k-1}\wedge \omega^{n-k},T-T'\rangle|.
\end{eqnarray*}

Since $u-u'\geq 0$ and $\pm\ddc\varphi_\epsilon$ are bounded by a constant times $\epsilon^{-2}\omega$, the second term in the last sum is bounded by a constant times $\epsilon^{-2}|\langle \mu_{k-1},u-u'\rangle|$. 
Applying the  Chern-Levine-Nirenberg inequality \cite[Ch.III (3.3)]{Demailly}  to the last term in the above sum, we see that this term is bounded by a constant times $\epsilon^\alpha$. This together with the induction hypothesis yields
$$|\langle \mu_k,u-u'\rangle|\lesssim \epsilon^{-2}|\langle \mu_{k-1},u-u'\rangle |+\epsilon^\alpha\lesssim \epsilon^{-2}\|u-u'\|_{L^1}^{\beta_{k-1}}+\epsilon^\alpha$$
for some constant $\beta_{k-1}>0$. 
Choosing  $\epsilon$ equal to a fixed constant small enough times $\|u-u'\|_{L^1}^{\beta_{k-1}\over 2+\alpha},$ the proof is  thereby complete. 

Note that we can show in the same way that a wedge-product of positive closed currents (of arbitrary bidegree) with H\"older continuous super-potential admits  a H\"older continuous super-potential.
\endproof

The following result, together with Theorem \ref{th_main}, implies Corollary \ref{main_cor}.

\begin{proposition} \label{prop_example}
Let $\mu$ be a positive measure on $X$. Assume that locally we can write $\mu=\ddc U+\partial V+\dbar W$ 
with H\"older continuous forms $U,V,W$ of bidegrees $(n-1,n-1)$, $(n-1,n)$ and $(n,n-1)$ respectively. Then $\mu$ admits a H\"older continuous super-potential. Moreover, the hypothesis on $\mu$ is satisfied when $\mu$ belongs to the Sobolev space $W^{2n/p-2+\epsilon,p}(X)$ or to the Besov space $B^{\epsilon-2}_{\infty,\infty}(X)$  for some $\epsilon>0$ and $p>1$. 
\end{proposition}
\proof
Consider a coordinate ball $\B$ in $X$ and $\chi$ a smooth positive function with compact support in $\B$.
We can assume that $\mu=\ddc U+\partial V+\dbar W$ as above on $\B$. Using a partition of unity, it is enough to show that $\mu':=\chi (\ddc U+\partial V+\dbar W)$ has a H\"older continuous super-potential. 

For $0<\epsilon\ll 1$, using the standard convolution, we can write $U=U_\epsilon+(U-U_\epsilon)$ with $\|U_\epsilon\|_{\Cc^2}\lesssim \epsilon^{-2}$ and $\|U-U_\epsilon\|_\infty\lesssim \epsilon^\alpha$ for some constant $\alpha>0$.  We obtain in the  same way  the  regularizing forms $V_\epsilon$ and  $W_\epsilon$
of  $V$ and  $W$ respectively. Using these  estimates,
 we have
\begin{eqnarray*}
|\langle \mu',u-u'\rangle| &\leq& |\langle \ddc U,\chi(u-u')\rangle| + |\langle \partial V,\chi(u-u')\rangle| + |\langle \dbar W,\chi(u-u')\rangle| \\
 &\leq& 
|\langle \ddc U_\epsilon, \chi(u-u')\rangle| +|\langle U-U_\epsilon, \ddc(\chi(u-u'))\rangle|\\
&& + |\langle \partial V_\epsilon, \chi(u-u')\rangle| +|\langle V-V_\epsilon, \partial(\chi(u-u'))\rangle|\\
&& + |\langle \dbar W_\epsilon, \chi(u-u')\rangle| +|\langle W-W_\epsilon, \dbar(\chi(u-u'))\rangle|.
       \end{eqnarray*}

Recall that $\|d u\|_{L^1}$ and $ \|d u'\|_{L^1}$ are bounded by a constant.  These properties are consequences of classical properties of psh functions. We can also obtain them
using   the estimates on $K$ and $\nabla K$ given in Proposition \ref{prop_norm_lp}.
Consequently,  $\ddc(\chi(u-u'))$, $ \partial(\chi(u-u'))$ and $\dbar(\chi(u-u'))$ have bounded mass.  This  discussion, combined with the  above estimates on $|\langle \mu',u-u'\rangle|$
and  the  above  mentioned properties of $U_\epsilon,$ $V_\epsilon,$ $W_\epsilon,$ implies that 
$$
 |\langle \mu',u-u'\rangle|\lesssim \epsilon^{-2} \|u-u'\|_{L^1} + \epsilon^\alpha.
$$
To obtain the first assertion in the proposition it is enough to take  $\epsilon:= \|u-u'\|_{L^1}^{1\over 2+\alpha}$. 

For the second assertion, locally on a small ball $\B$ in $X$ with holomorphic coordinates $z$, if $u$ is a solution of the Laplacian equation $\Delta u=\mu$ and if $U$ is a suitable constant times $u(\ddc\|z\|^2)^{n-1}$, then  $\ddc U=\mu$. When $\mu$ belongs to the Sobolev space $W^{2n/p-2+\epsilon,p}(X)$, the function $u$ is in $W^{2n/p+\epsilon,p}(\B)$, see \cite[p.198]{Brezis} for $W^{k,p}$ with $k\in\N$ and \cite[p.186 and p.230]{Triebel} for the interpolation and the duality which allow to consider the case of $W^{k,p}$ with $k\in\R$. By Sobolev's embedding theorem, $u$ is H\"older continuous, see e.g. \cite[p.168]{Brezis}. When $\mu$ is in the Besov space $B^{\epsilon-2}_{\infty,\infty}(X)$, the solution $u$ is in $B^{\epsilon}_{\infty,\infty}$ which is the H\"older space $\Cc^\epsilon$.  The result follows.
\endproof

The following result and Theorem \ref{th_main} imply the main result of Hiep in \cite{Hiep}.

\begin{proposition}
Let $\mu$ be a positive measure on $X$. Assume that there are constants $c>0$ and $\alpha>0$ such that if $B$ is a ball of radius $r$ in $X$ we have $\mu(B)\leq cr^{2n-2+\alpha}$. Then $\mu$ has a H\"older continuous super-potential.
\end{proposition}
\proof
Let $T,T',u,u'$, the super-potential $\Uc$, the kernel $K,$  the  function $v$ and the constant $m$ be given in Proposition \ref{prop_norm_lp}  above. We have 
\begin{eqnarray*} \Uc(T)-\Uc(T') & = & \langle \mu,u-u'\rangle =  \langle \mu,v-m\rangle\\
&=&
\big\langle \mu,(\pi_1)_*(K\wedge\pi_2^*(T-T')) \big\rangle -m\\
& = & \big\langle T-T',(\pi_2)_*(K\wedge \pi_1^*(\mu))\big\rangle-m.
\end{eqnarray*}
Since we already have good estimates on $m$, 
it is enough to check that the form $\Phi:=(\pi_2)_*(K\wedge \pi_1^*(\mu))$ is H\"older continuous.

Let $\B$ be a coordinate ball in $X$ that we identify with the unit ball in $\C^n$. We will show that $\Phi$ is H\"older continuous near the origin $0\in\C^n$.  Let $\chi$ be a smooth function with compact support in $\B\times \B$ and equal to 1 near the origin. We have 
$$\Phi= (\pi_2)_*(\chi K\wedge \pi_1^*(\mu))+ (\pi_2)_*((1-\chi)K\wedge \pi_1^*(\mu)).$$
Since the form $(1-\chi)K$ is smooth near $X\times \{0\}$, the last expression in the above identity defines a smooth form near 0. It remains to show that $\Psi:=(\pi_2)_*(\chi K\wedge \pi_1^*(\mu))$ is H\"older continuous near 0.

Observe that the coefficients of $\Psi$ have the form
$$\Theta(x):=\int_{y\in \B} H(x,y)d\mu(y)$$
where $H$ is a coefficient of $\chi K$. We use now the estimates on $K$ and $\nabla K$ given in Proposition \ref{prop_norm_lp}. Consider two points $x$ and $x'$ in $\B$ near 0 and denote by $D$ the ball of center $x$ and of radius $\rho:=2\|x-x'\|^{1/(2n)}$. We have 
$$\Theta(x)-\Theta(x')=\int_DH(x,y)d\mu(y) -\int_D H(x',y)d\mu(y) + \int_{\B\setminus D} (H(x,y)-H(x',y))d\mu(y).$$
The estimate on $K$ and the hypothesis on $\mu$ imply that the first two terms are bounded by a constant times $|\log\rho |\rho^{\alpha}$. The estimate on $\nabla K$ implies that 
$$|H(x,y)-H(x',y)|\lesssim \|x-x'\| \rho^{1-2n} \quad \mbox{for } y\not\in D.$$ 
Therefore, the last integral  in the above identity is bounded by a constant times $\|x-x'\|^{1/(2n)}$. We deduce that  $\Theta$ is a H\"older continuous function. The proposition follows. 
\endproof

\begin{proposition} \label{prop_holder_moderate}
Let $\mu$ be a positive measure on $X$ with a H\"older continuous super-potential. Then $\mu$ is moderate. 
\end{proposition}
\proof
Let  $u$ be  an $\omega$-psh function with $\int_X u\omega^n=0$  and set $u_M:=\max(u,-M)$ for $M>1$ large enough. 
By  Lemma \ref{lemma_moderate_Lp} applied to the measure $\omega^n$, there is  a constant  $\alpha>0$ independent of $u$ and  $M$ such that
$$
\Big|  \int_X u_M\omega^n\Big|= \Big|  \int_X (u_M-u)\omega^n\Big|\lesssim e^{-\alpha M}.
$$

Denote by $\beta$ a H\"older exponent of $\Uc$ with respect to the distance  $\dist_{L^1}$ on $\Cc$ with $0<\beta\leq 1$. Since $M$ is large enough, the estimate above shows that 
$\|u-u_M\|_{L^1}$ is small. We deduce from Lemma \ref{lemma_elem}  and the inequality $u_M-u >1$ on  $\{u<-M-1\}$ that
$$\mu\{u <-M-1\} \leq   \Big|\int (u-u_M)d\mu          \Big| \lesssim \| u-u_M\|^\beta_{L^1}\lesssim e^{-\alpha \beta M}.$$
Thus,  $\mu$ is  moderate.
\endproof

The following corollary gives us a large family  of measures with H\"older continuous super-potential. It shows that the restriction of such a measure to a Borel set 
gives also a measure with H\"older continuous super-potential. Note that by definition the
set of measures with H\"older continuous super-potential is a convex cone. We then deduce from Theorem \ref{th_main} analogous properties for Monge-Amp\`ere measures with H\"older continuous potential which have been obtained in \cite{DemaillyDinew}. 

\begin{corollary} \label{cor_example}
Let $\mu$ be a positive measure with a H\"older continuous super-potential. If $f$ is a positive function in $L^p(\mu)$ with $p>1$, then $f\mu$ admits a H\"older continuous super-potential.
\end{corollary}
\proof
Let $q\geq 1$ be the real number such that $p^{-1}+q^{-1}=1$.  Let $\Vc$ be the super-potential of $f\mu$. 
Consider two currents $T,T'$ in $\Cc$ and their $\omega$-potentials $u,u'$.
We have
$$|\Vc(T)-\Vc(T')|=|\langle f\mu,u-u'\rangle|\leq \|f\|_{L^p(\mu)}\|u-u'\|_{L^q(\mu)}.$$
By Proposition \ref{prop_holder_moderate} and Lemma \ref{lemma_moderate_Lp}, the last expression is bounded by a constant times $\|u-u'\|_{L^1(\mu)}^{1/(2q)}$. 
Lemma \ref{lemma_elem} implies the result.
\endproof

\noindent
{\bf End of the proof of Theorem \ref{th_main}.} The necessary condition follows from Proposition \ref{prop_necessary}. 
To prove   the  sufficient part 
assume that the super-potential $\Uc$ of $\mu$ is H\"older continuous with H\"older exponent $0<\beta\leq 1$ with respect to the distance $\dist_{L^1}$ on $\Cc$.
Propositions \ref{prop_holder_moderate} and \ref{prop_moderate_K} imply that for every $p>0$ there is a constant $c>0$ such that 
$$\mu(A)\leq c\capacity_\BTK(A)^p \quad \mbox{for any compact set } A.$$
In other  words, $\mu$ satisfies the condition $\cal H(\infty)$ in the sense of \cite{EGZ} (see also Definition 2.6 in \cite{DemaillyDinew})
which  is  stronger than condition (A) in \cite{Kolodziej1}. The latter work ensures the existence  of a bounded $\omega$-psh solution
$u$ to  $(\ddc u+\omega)^n=\mu$ with the  normalization $\inf_X u=1$ (for the reader's convenience we use here the notation $u$ instead of $\varphi$ as in the reference \cite{DemaillyDinew} that we use now).
 
For $\delta>0$ small enough,  let $\rho_\delta u$ denote  the  regularization of $u$  defined in  \cite[Sec. 2]{DemaillyDinew}, see also \cite{Demailly2,Demailly3}. 
It satisfies 
$$\|\rho_\delta u-u\|_{L^1} \leq c \delta^2 \quad \mbox{and} \quad 
\ddc \rho_\delta u\geq -c\omega$$
for some constant $c\geq 1$. 
Lemma \ref{lemma_elem} applied to the $\omega$-psh functions $c^{-1}\rho_\delta u$ and $c^{-1}u$ implies that
\begin{equation} \label{eq_crucial_estimate}
\int_X |\rho_\delta u-u|d\mu\lesssim \|\rho_\delta u-u\|_{L^1}^\beta\lesssim  \delta^{2\beta}.
\end{equation}
We  are now able  to apply Proposition 3.3 in \cite{DemaillyDinew} which shows that $u$ is  H\"older continuous.

Now, we follow  along the same lines as those given in the proof of Theorem A in \cite{DemaillyDinew}
in order to  get  an explicit H\"older exponent of $u$
in terms of $n$ and $\beta$. It is enough to make the  following changes:
\begin{itemize}
\item[$\bullet$]  After (3.1) in  \cite{DemaillyDinew}:  let $q:=1$,  $0<\alpha_1<{2\beta\over n+1}$ and choose 
$\epsilon>0,$   $\alpha,$  $\alpha_0$ such that 
$ \alpha_1<\alpha<\alpha_0<2\beta-\alpha_0 (n+\epsilon).$ Take $f\equiv 1$ and set $g:=0$ on $E_0$  and $g:=c$ elsewhere with a constant  $c\geq 1$  such that $\|g\mu\|=\|\mu\|$.
\item[$\bullet$] As in  \cite{DinewZhang, Kolodziej2}, we  solve for continuous $\omega$-psh function $v$ 
$$(\ddc v+\omega)^n=g\mu \quad \mbox{with} \quad \max(u-v)=\max(v-u). $$
\item[$\bullet$]  (3.2) in  \cite{DemaillyDinew} becomes now (\ref{eq_crucial_estimate}) which implies 
$\int_{E_0} d\mu \leq  c_2 \delta^{2\beta-\alpha_0}$ with $c_2>0$. The last relation replaces 
the inequality $\int_{E_0} f\omega^n\leq  c_2 \delta^{2-\alpha_0\over q}$ in \cite{DemaillyDinew}.
\end{itemize}
We obtain that $u$ is H\"older continuous with  exponent $\alpha_1$ for all   $0<\alpha_1<{2\beta\over n+1}.$
\hfill $\square$

\noindent
T.-C. Dinh, UPMC Univ Paris 06, UMR 7586, Institut de
Math{\'e}matiques de Jussieu, 4 place Jussieu, F-75005 Paris, France.\\ 
{\tt dinh@math.jussieu.fr}, {\tt http://www.math.jussieu.fr/$\sim$dinh}

\medskip

\noindent
V.-A.  Nguy{\^e}n,
Math{\'e}matique-B{\^a}timent 425, UMR 8628, Universit{\'e} Paris-Sud,
91405 Orsay, France.\\
  {\tt VietAnh.Nguyen@math.u-psud.fr}, {\tt http://www.math.u-psud.fr/$\sim$vietanh}

\end{document}